\documentclass{amsart}
\usepackage{amsmath}
\usepackage{graphicx}
\usepackage{amsfonts}
\usepackage{amssymb}%
\providecommand{\U}[1]{\protect\rule{.1in}{.1in}}
\providecommand{\U}[1]{\protect\rule{.1in}{.1in}}
\providecommand{\U}[1]{\protect\rule{.1in}{.1in}}
\providecommand{\U}[1]{\protect\rule{.1in}{.1in}}
\providecommand{\U}[1]{\protect\rule{.1in}{.1in}}
\newtheorem{theorem}{Theorem}

\newtheorem{corollary}[theorem]{Corollary}

\newtheorem{definition}[theorem]{Definition}

\newtheorem{proposition}[theorem]{Proposition}
\newtheorem{remark}[theorem]{Remark}

\begin{document}
\title{ Index of Graded Filiform and Quasi Filiform Lie Algebras.}
\author{Hadjer  Adimi}
\address{H. Adimi, Universit\'e de Bordj Bou Arreridj, Algeria}
\email{adimiha2001@yahoo.fr}
\author{Abdenacer Makhlouf}
\address{A. Makhlouf, Laboratoire du Math\'{e}matiques, Informatiques  et applications, Universit\'{e} de
Haute-Alsace, 4 rue des Fr\`{e}res Lumi\`{e}re F-68093 Mulhouse cedex, France.}
\email{Abdenacer.Makhlouf@uha.fr}

\begin{abstract}
The \ filiform and the quasi-filiform Lie algebras form a special class of
nilpotent Lie algebras. The aim of this paper is to compute the index and
provide regular vectors of this two class of nilpotent Lie algebras. we
consider the graded filiform Lie algebras $L_{n}$, $Q_{n}$, the $n$-dimensional filiform Lie algebras for $n$ $<8$, also the graded
quasi-filiform Lie algebras and finally a Lie algebras whose nilradical is
$Q_{2n}$.

\end{abstract}

\maketitle

\section{Introduction}
The index of Lie algebra has applications to invariant theory and  are of interest in deformation and quantum group theory. A Lie algebra is said to be Frobenius if the index is $0$ which is equivalent to say that there is a functional in the dual such that the bilinear form $B_F$, defined by $B_F(x,y)=F([x,y])$, is nondegenerate. Frobenius algebras were first studied by Ooms in \cite{Al}. He proved that the universal enveloping algebra of the Lie algebra is primitive, that is it admits a faithful simple module, provided that the Lie algebra is Frobenius and that the converse  holds when the Lie algebra is algebraic. Most of the studies of index concerned  simple Lie algebras or  their subalgebras. They have been considered by many authors \cite{De-Ki,El1,El2,El3,El4,Al2008,Pa-Ru1,Pa-Ru2}. Notice that simple Lie algebra can never be Frobenius but many subalgebras are. In this paper we focus on the computation of the index for nilpotent Lie algebras, mainly the class of filiform and quasi-filiform Lie algebras. 

The filiform Lie algebras were introduced by M. Vergne (see \cite{vergne}),
she classified them up to dimension 6 and also characterized the graded
filiform Lie algebras. In the last decades several papers were dedicated to
classification of filiform Lie algebras of larger dimensions.  In particular
$L_{n}$ plays an important role in the study of filiform and nilpotent Lie
algebras.  It is known that any $n$-dimensional filiform Lie algebra may be
obtained by deformation of the one of the filiform Lie algebras $L_{n}$. The
classification of naturally graded quasi-filiform Lie algebras is known.  They
have the characteristic sequence $\left(  n-2,1,1\right)  $ where $n$ is the
dimension of the algebra. The aim of this paper is to  give an extended version of
our paper \cite{AdimiMakhloufLie} and to focus on filiform Lie algebras and quasi-filiform Lie algebras. We
compute the index and provide the regular vectors of $n$-dimensional filiform
Lie algebras for $n<8$ and quasi-filiform Lie algebras. In the first Section,
we summarize the index theory of Lie algebras. Then, in Section 2, we review
the nilpotent and filiform Lie algebras theories. Section 3, is dedicated to
the two graded filiform Lie algebras $L_{n}$ and $Q_{n}$. In
Section 4, we consider the classification up to dimension 8 and compute for
each filiform Lie algebra its index and the set of all regular vectors. In
section 5 we compute the index of graded quasi-filiform Lie algebras, and provide corresponding 
regular vectors. In the last section we compute the index of Lie
algebras whose nilradical is $Q_{2n}$.

\section{Lie algebras Index}
Throughout this paper $\mathbb{K}$ is an algebraically closed field of
characteristic 0.
In this Section, we summarize the index theory of Lie algebras. Let
$\mathcal{G}$ be an $n$-dimensional  Lie algebra. Let $x\in\mathcal{G}$, we denote by $adx$ the
endomorphism of $\mathcal{G}$ defined by $adx\left(  y\right)  =\left[
x,y\right]  $ for all  $y\in\mathcal{G}.$

\begin{definition}
A Lie algebras $\mathcal{G}$ over $\mathbb{K}$ is a pair consisting of a
vector space $\mathbb{V}=\mathcal{G}$ and a skew-symmetric bilinear map
$[~,~]\ :\ \mathcal{G}\ \times\mathcal{G}\ \rightarrow\mathcal{G}\ \ \left(
x,y\right)  \rightarrow\left[  x,y\right]  \ $ satisfying the Jacobi identity
\[
\left[  x,\left[  y,z\right]  \right]  +\left[  y,\left[  z,x\right]  \right]
+\left[  z,\left[  x,y\right]  \right]  =0\quad\mathit{\ \ \ }\forall
x,y,z\in\mathcal{G}.
\]

\end{definition}

Let $\mathbb{V}$ be a finite-dimensional vector space over $\mathbb{K}$
provided with the Zariski topology, $\mathcal{G}$ be a Lie algebra and
$\mathcal{G}^{\ast}$ its dual. Then $\mathcal{G}$ actes on $\mathcal{G}^{\ast
}$ as follows:
\begin{align*}
\mathcal{G\times G}^{\ast}  &  \rightarrow\mathcal{G}^{\ast}\\
\left(  x,f\right)   &  \mapsto x\cdot f
\end{align*}
where $\mathit{\forall\in G:}\left(  x\cdot f\right)  \left(  y\right)
\mathit{=f}\left(  \left[  x,y\right]  \right)  \mathit{.}$

Let $f\in\mathcal{G}^{\ast}$ and $\Phi_{f}$ be a skew-symmetric bilinear form
defined by
\begin{align*}
\Phi_{f}:\mathcal{G\times G}  &  \rightarrow\mathbb{K}\\
\left(  x,y\right)   &  \mapsto\Phi_{f}\left(  x,y\right)  =f\left(  \left[
x,y\right]  \right).
\end{align*}
We denote the kernel of the map $\Phi_{f}$ by $\mathcal{G}^{f}$,
\begin{equation}
\mathcal{G}^{f}=\{x\in\mathcal{G}:f([x,y])=0\ \ \forall y\in\mathcal{G}\}.
\end{equation}

\begin{definition}
The index of a Lie algebra $\mathcal{G}$ is the integer
\[
\chi_{\mathcal{G}}=\inf\left\{  \dim\text{ }\mathcal{G}^{f};f\in
\mathcal{G}^{\ast}\right\}  .
\]

A linear functional $f\in\mathcal{G}^{*}$ is called regular if $\mathit{dim }
\ \mathcal{G}^{f }=\chi_{\mathcal{G}}$. The set of all regular linear
functionals is denoted by $\mathcal{G}_{r}^{*}.$
\end{definition}

\begin{remark}
The set $\mathcal{G}_{r}^{*}$ \ of all regular linear functionals is a
nonempty Zariski open set.
\end{remark}

Let $\{x_{1},\cdots,x_{n}\}$ be a basis of $\mathcal{G}$. We can express the
index using the matrix $([x_{i},x_{j}])_{1\leq i<j\leq n}$ as a matrix over
the ring $S(\mathcal{G})$, (see \cite{Dixmier}). We has the following
proposition.

\begin{proposition}
\label{Dixmier} The index of an $n$-dimensional Lie algebra $\mathcal{G}$ is
the integer
\[
\chi\left(  \mathcal{G}\right)  =n-Rank_{R(\mathcal{G})}\left(  \left[
x_{i},x_{j}\right]  \right)  _{1\leq i<j\leq n}%
\]
where $R(\mathcal{G})$ is the quotient field of the symmetric algebra
$S(\mathcal{G})$.
\end{proposition}

\begin{remark}
The index of an $n$-dimensional abelian Lie algebra is $n.$
\end{remark}

\begin{proposition}
\label{CentralExtension}Let $\mathcal{G}_{0}$ be a Lie algebra, and
$\mathcal{G}$ be a central extension of $\mathcal{G}_{0}$ by a $1-$dimensional
Lie algebra $\mathcal{L=}c\mathbb{C}$, then $\chi\left(  \mathcal{G}\right)
=\chi\left(  \mathcal{G}_{0}\right)  +1.$ Moreover $f$ is a regular vector of
$\mathcal{G}$ then $f=g+\rho c^{\ast}$ where $g$ is a regular vector of
$\mathcal{G}_{0}$ and $\rho\in\mathbb{C}.$
\end{proposition}

\begin{proof}
Indeed, we have
$$\left\{
\begin{array}
[c]{c}%
\left[  x,c\right]  =0\ \ \ \forall x\in\mathcal{G}_{0},\\
\left[  c,c\right]  =0.\
\end{array}
\right.  $$
Then the matrix associated to $\mathcal{G}$ is of the form
$M=\left(
\begin{array}
[c]{cc}%
M_{\mathcal{G}_{0}} & 0\\
0 & 0
\end{array}
\right)  .$ It follows that  $Rank\left(  \mathcal{G}\right)  =Rank\left(
\mathcal{G}_{0}\right)  .$ Therefore
\[
\chi\left(  \mathcal{G}\right)  =\chi\left(  \mathcal{G}_{0}\right)  +1.
\]

Let $g$ be a regular vector of $\mathcal{G}_{0}$.  Then $\dim$ $\mathcal{G}%
_{0}^{g}=\chi\left(  \mathcal{G}_{0}\right)  $ and $f=g+\rho c^{\ast}$ is a
regular vector of $\mathcal{G}$.

We know that $\mathcal{G}^{f}=\left\{  x\in\mathcal{G},\ \ f\left(
\left[  x,y\right]  \right)  =0,\forall y\in\mathcal{G}\right\}  .$

We set
\[
x=x_{0}+\lambda c\text{ and }y=y_{0+}\mu c.
\]
Then 
\begin{eqnarray*}
f\left(  \left[  x,y\right]  \right) =g\left(  \left[  x,y\right]
\right)  +\rho c^{\ast}(\left[  x,y\right]  )= g\left(  \left[  x_{0},y_{0}\right]  \right). 
\end{eqnarray*}
We have
\[
g\left(  \left[  x_{0},y_{0}\right]  \right)  =0\text{ if }x_{0}\in
\mathcal{G}_{0},\forall y\in\mathcal{G}_{0}.
\]
Therefore $\mathcal{G}^{f}=\mathcal{G}_{0}^{g}+c\mathbb{C}$.
\end{proof}

\begin{remark}
\label{regularvector}In the sequel, we use the following procedure to compute  regular vectors. 
We recall that if $\dim\mathcal{G}^{f}=\chi\left(  \mathcal{G}\right)  $ then
$f$ is a regular vector of $\mathcal{G}$, where
$\chi\left(  \mathcal{G}\right)  =\min\left\{  \dim\mathcal{G}^{f}%
,f\in\mathcal{G}^{\ast}\right\}  $ and
$\mathcal{G}^{f}=\{x\in\mathcal{G}:f([x,y])=0\ \forall y\in\mathcal{G}\}.$

The equation
$f([x,y])=0$ implies $\sum_{i=1}^{n}\sum_{j=1}^{n}\sum_{s=1}^{n}{a_{i}
b_{j}p_{s}x_{s}^{\ast}}\left(  [{x_{i}},{x}_{j}]\right)  =0.$ 

It is equivalent to  

$\sum_{i=1}^{n}\sum_{j=1}^{n}\sum_{s=1}^{n}{a_{i}b_{j}p_{s}C}_{ij}^{s}=0,$ where $C_{ij}^{s}$ are the structure constants with respect to the basis $\{ x_i\}_i$. 
Then for all  $
j$, we have $\sum_{s=1}^{n}\sum_{i=1}^{n}{a_{i}p_{s}C}_{ij}^{s}=0.$
It leads to 
$$\left(
\sum_{s=1}^{n}{p_{s}C}_{ij}^{s}\right)  _{ij}\left(
\begin{array}
[c]{c}%
a_{1}\\
.\\
.\\
.\\
a_{n}%
\end{array}
\right)  =\left(
\begin{array}
[c]{c}%
0\\
.\\
.\\
.\\
0
\end{array}
\right)  $$

We denote by $M=\left(  \sum_{s=1}^{n}{p_{s}C}_{ij}^{s}\right)  _{ij}$ and assume ${C}_{ij}^{s}=-{C}_{ji}^{s}$

We search the minors of order $n-\chi\left(  \mathcal{G}\right)  $ of non-zero
determinant of the matrix $M$.

The matrix $M=\left(  \sum_{s=1}^{n}{p_{s}C}_{ij}^{s}\right)  _{ij}$ is the
same matrix as the multiplication table in which we replace ${x_{s}}$ by ${p_{s}}$.
\end{remark}

\begin{definition}
A Lie algebra $\mathcal{G}$ over an algebraically closed field of
characteristic $0$ is said to be Frobenius if there exists a linear form
$f\in\mathcal{G}^{\ast}$ such that the bilinear form $\Phi_{f}$ on
$\mathcal{G}$ is nondegenerate.
\end{definition}

In \cite{El3}, the author described all the Frobenius algebraic Lie algebras
$\mathcal{G}=R+N$ whose nilpotent radical $N$ is abelian in the following two
cases: the reductive Levi subalgebra $R$ acts on $N$ irreducibly and  $R$ is
simple. He classified all the algebraic Frobenius algebras up to dimension 6.
See also \cite{Al} and \cite{Al2008} for further computations.

We discuss now the evolution by deformation of the index of a Lie algebra.
About deformation theory we refer to \cite{Ge,NR,MakDeform}. Let
$\mathbb{V}$ be a $\mathbb{K}$-vector space and $\mathcal{G}_{0}=(\mathbb{V},
[~, ~]_{0})$ be a Lie algebra. Let $\mathbb{K }[[t]]$ be the power series ring
in one variable $t$ and coefficients in $\mathbb{K }$ and $\mathbb{V}[[t]]$ be
the set of formal power series whose coefficients are elements of $\mathbb{V}%
$. A \emph{formal Lie deformation} of $\mathcal{G}_{0}$ is given by the
$\mathbb{K}[[t]]$-bilinear map $[~, ~]_{t}: \mathbb{V}[[t]] \times
\mathbb{V}[[t]] \rightarrow\mathbb{V}[[t]] $ of the form $[~,~]_{t}
=\sum_{i\geq0}[~,~]_{i}t^{i}, $ where each $[~,~]_{i}$ is a $\mathbb{K}%
$-bilinear map $[~,~]_{i}: \mathbb{V}\times\mathbb{V} \rightarrow\mathbb{V}$,
satisfying the skew-symmetry and the Jacobi identity.

\begin{proposition}
\label{propDeformation} The index of a Lie algebra decreases by one parameter formal deformation.
\end{proposition}

\begin{proof}
The rank of the matrix $\left(  \left[  x_{i},x_{j}\right]  \right)  _{ij}$
increases by deformation, consequently the index decreases.
\end{proof}

\section{Nilpotent and Filiform Lie algebras}

 In this Section, we review the theory of nilpotent and
filiform Lie algebras.
Let $\mathcal{G}$ be a Lie algebra. We set $\mathcal{C}^{0}\mathcal{G}%
=\mathcal{G}$ and $\mathcal{C}^{k} \mathcal{G}=[\mathcal{C}^{k-1}%
\mathcal{G},\mathcal{G}]$, for $k>0$. A Lie algebra $\mathcal{G}$ is said to
be nilpotent if there exists an integer $p$ such that $\mathcal{C}%
^{p}\mathcal{G}=0$. The smallest $p$ such that $\mathcal{C}^{p}\mathcal{G}=0$
is called the nilindex of $\mathcal{G}$. Then, a nilpotent Lie algebra has a
natural filtration given by the central descending sequence:
\[
\mathcal{G}=\mathcal{C}^{0}\mathcal{G}\supseteq\mathcal{C}^{1}\mathcal{G}%
\supseteq\cdots\mathcal{C}^{p-1} \mathcal{G}\supseteq\mathcal{C}^{p}
\mathcal{G}=0.
\]
We have the following characterization of nilpotent Lie algebras (Engel's Theorem).

\begin{theorem}
A Lie algebra $\mathcal{G}$ is nilpotent if and only if the operator $adx$ is
nilpotent for all $x$ in $\mathcal{G}$.
\end{theorem}

In the study of nilpotent Lie algebras,  filiform Lie algebras, which were introduced by M. Vergne, play an
important role.  An $n$-dimensional
nilpotent Lie algebra is called \emph{filiform} if its nilindex $p$ equals $n-1$. The
filiform Lie algebras are the nilpotent algebras with the largest nilindex. If
$\mathcal{G}$ is an $n$-dimensional filiform Lie algebra, then we have $dim
\mathcal{C}^{i} \mathcal{G}=n-i \quad\text{for }\ \ 2\leq i\leq n.$

Another characterization of filiform Lie algebras uses characteristic
sequences $c(\mathcal{G})=sup\{c(x):x\in\mathcal{G}\setminus\lbrack
\mathcal{G},\mathcal{G}]\},$ where $c(x)$ is the sequence, in decreasing
order, of dimensions of characteristic subspaces of the nilpotent operator
$adx$.

\begin{definition}
An $n$-dimensional nilpotent Lie algebra is filiform if its characteristic
sequence is of the form $c\left(  \mathcal{G}\right)  =\left(  n-1,1\right)
.$
\end{definition}
\section{Index of Graded filiform Lie algebras}

The classification of filiform Lie algebras was given by Vergne (\cite{vergne}%
) up to  dimension 6 and then extended to dimension 11 by several authors (see
\cite{AncocheaGozeNilpo,CarlesNilpo,GozeKhakimjanov,Seeley,Gomez1}).

In the general case there is two classes $L_{n}$ and $Q_{n}$ of filiform Lie
algebras which plays an important role in the study of the algebraic varieties
of filiform and more generally nilpotent Lie algebras.

Let $\{x_{1},\cdots,x_{n}\}$ be a basis of the $\mathbb{K}$ vector space
$L_{n}$, the Lie algebra structure of $L_{n}$ is defined by the following
non-trivial brackets : $\left[  x_{1},x_{i}\right]  =x_{i+1}\quad
i=2,...,n-1.$

Let $\{x_{1},\cdots,x_{n=2k}\}$ be a basis of the $\mathbb{K}$ vector space
$Q_{n}$, the Lie algebra structure of $Q_{n}$ is defined by the following
non-trivial brackets.
\begin{align*}
Q_{n}:\left[  x_{1},x_{i}\right]   &  =x_{i+1}\quad i=2,...,n-1,\\
\quad\quad\left[  x_{i},x_{n-i+1}\right]   &  =\left(  -1\right)  ^{i+1}%
x_{n}\quad i=2,...,k\quad\text{ where }n=2k.
\end{align*}

The classification of $n$-dimensional graded filiform Lie algebras yields to
two isomorphic classes $L_{n}$ and $Q_{n}$ when $n$ is odd and to only the Lie
algebra $L_{n}$ when $n$ is even.

It turns out that any filiform Lie algebra is isomorphic to a Lie algebra
obtained as a deformation of a Lie algebra $L_{n}$.

\subsection{Index of Filiform Lie algebras}

We aim to compute the index of $L_{n}$ and $Q_{n}$ and regular vectors.

\paragraph{Index of $L_{n}$ :}\

Let $\left\{  x_{1},x_{2},...,x_{n}\right\}  $ be a fixed basis of the vector
space $\mathbb{V}=L_{n}$ and $\left\{  x_{1}^{\ast},.....,x_{n}^{\ast
}\right\}  $ be a basis of the dual space. Set $f=\sum_{i\geq1}{p_{i}%
x_{i}^{\ast}}\in\mathbb{V}^{\ast}.$

\begin{proposition}
For $n\geq3$, the index of the $n$-dimensional filiform Lie algebra $L_{n}$ is
$\chi\left(  L_{n}\right)  =n-2.$ The regular vectors of $L_{n}$ are of the
form $f={p_{1}\ x_{1}^{\ast}+}$ ${p_{2}\ x_{2}^{\ast}+p_{s}\ x_{s}^{\ast}}$
where $s\in\{3,...,n\}$ and ${p_{s}\neq0.}$
\end{proposition}

\begin{proof}
Since the corresponding matrix to the Lie algebra $L_{n}$ is of the form
\[
\left(
\begin{array}
[c]{ccccc}%
0 & x_{3} & \cdots & x_{n} & 0\\
-x_{3} & 0 & \cdots & 0 & 0\\
\vdots & \vdots & \vdots & \vdots & \vdots\\
-x_{n} & 0 & \cdots & 0 & 0\\
0 & 0 & \cdots & 0 & 0
\end{array}
\right)
\]

and its rank is 2, then $\chi(L_{n})=n-2.$ The second assertion is obtained by
a direct calculation:

We set  $x=\sum_{i=1}^{n}{\ a_{i}x_{i},}$ $y=\sum_{j=1}^{n}{\ b_{j}x_{j},}$
$f=\sum_{s=1}^{n}p{_{s}x_{s}^{\ast}}$
and $\mathcal{G}^{f}=\{x\in\mathcal{G}:f([x,y])=0\ \ \forall y\in
\mathcal{G}\}.$

Then $f([x,y])=0$ implies
$\sum_{i=1}^{n}\sum_{j=1}^{n}\sum_{s=1}^{n}{a_{i}b_{j}p_{s}x_{s}^{\ast}
}\left(  [{x_{i}},{x}_{j}]\right)  =0.$
 It is equivalent to
$$\sum_{s=1}^{n}\sum_{j=2}^{n-1}{a_{1}b_{j}p_{s}x_{s}^{\ast}}\left(  [{x_{1}%
},{x}_{j}]\right)  -\sum_{i=2}^{n-1}{a_{i}b_{1}p_{s}x_{s}^{\ast}}\left(
[{x_{1}},{x}_{j}]\right)  =0.$$
 Then we obtain $\sum_{s=1}^{n}\sum_{i=2}^{n-1}\left(  {a_{1}b_{i}-a_{i}b_{1}}\right)
{p_{s}x_{s}^{\ast}}\left(  {x_{i+1}}\right)  =0.$ The equation 
$\sum_{i=2}^{n-1}\left(  {a_{1}b_{i}-a_{i}b_{1}}\right)  {p_{i+1}}=0$ should hold 
for all ${b_{i}}$. 
It leads to the following  system
$$
\begin{cases}
{a_{1}p_{i+1}=0, \ \ 2\leq i\leq n-1},\\
\sum_{i=2}^{n-1}{a_{i}p_{i+1}=0}.
\end{cases}
  $$

Therefore, 
one of the $p_{i}$ satisfies $p_{i}\neq0$ where i$\in\{3,...,n\}.$
\end{proof}

\paragraph{Index of $Q_{n}$ :}

\begin{proposition}
For $n=2k$ and $k\geq2$, the index of the $n$-dimensional filiform Lie algebra
$Q_{n}$ is $\chi\left(  Q_{n}\right)  =2.$ The regular vectors of $Q_{n}$ are
of the form $f=\sum_{i=1}^{n}{\ p_{i}x_{i}^{\ast}}$ with $p_{n}\neq0$.
\end{proposition}

\begin{proof}
Since the corresponding matrix to the Lie algebra $Q_{n}$ is of the form

$$\left(
\begin{array}
[c]{ccccccc}%
0 & x_{3} & x_{4} & \cdots & x_{n-1} & x_{n} & 0\\
-x_{3} & 0 & 0 & \cdots & 0 & -x_{n} & 0\\
-x_{4} & 0 & 0 & \cdots & x_{n} & 0 & 0\\
\vdots & \vdots & \vdots & \vdots & \vdots & \vdots & \vdots\\
-x_{n-1} & 0 & -x_{n} & \cdots & 0 & 0 & 0\\
-x_{n} & x_{n} & 0 & \cdots & 0 & 0 & 0\\
0 & 0 & 0 & 0 & 0 & 0 & 0
\end{array}
\right)  $$

and its rank is n-2, then $\chi(Q_{n})=2.$ The second assertion is obtained by
the following  calculation.

Let $\left\{  x_{1},x_{2},.....,x_{n}\right\}  $ be a fixed of  basis of $Q_{n}$,
$x=\sum_{i=1}^{n}{\ a_{i}x_{i}}$, $y=\sum_{j=1}^{n}{\ b_{j}x_{j},}$
$f=\sum_{s=1}^{n}{\ p_{s}x_{s}^{\ast}}$
 and $\mathcal{G}^{f}=\{x\in\mathcal{G}:f([x,y])=0\ \ \forall y\in\mathcal{G}\}. $

The equation  $f([x,y])=0$ may be written as 
$\sum_{i=1}^{n}\sum_{j=1}^{n}\sum_{s=1}^{n}{a_{i}b_{j}p_{s}x_{s}^{\ast}%
}\left(  [{x_{i}},{x}_{j}]\right)  =0.$ It is equivalent to 
$\sum_{i=2}^{n-1}\left(  {a_{1}b}_{i}{-a_{i}b_{1}}\right)  {p_{i+1}+}%
\sum_{i=2}^{n-1}\left(  -1\right)  ^{i+1}\left(  {a_{1}b_{n-i+1}%
-a_{n-i+1}b_{i}}\right)  {p_{n}}=0$. Then 
$$
\begin{cases}
\sum_{i=2}^{n-1}\left(  {a_{1}p_{i+1}}\right)  {b_{i}=0},\\
-{b_{1}}\sum_{i=2}^{n-1}{a_{i}p_{i+1}=0},\\
\sum_{i=2}^{n-1}\left(  -1\right)  ^{i+1}\left(  {a_{1}p_{n}}\right)
{b_{n-i+1}=0},\\
-\sum_{i=2}^{n-1}\left(  -1\right)  ^{i+1}\left(  {a_{n-i+1}p_{n}}\right)
{b_{i}=0}.
\end{cases}
  $$

Canceling  the first and the last columns and the corresponding lines, leads to the following minor

$$\left(
\begin{array}
[c]{ccccc}%
0 & 0 & \cdots & 0 & -x_{n}\\
0 & 0 & \cdots & x_{n} & 0\\
\vdots & \vdots & \vdots & \vdots & \vdots\\
0 & -x_{n} & \cdots & 0 & 0\\
x_{n} & 0 & \cdots & 0 & 0
\end{array}
\right)  $$

Hence, we obtain   $f=\sum_{i=1}^{n}{\ p_{i}x_{i}^{\ast},}$ with $p_{n}\neq0$.
\end{proof}

Using Proposition \ref{propDeformation}, we obtain the following result.

\begin{corollary}
The index of a filiform Lie algebra is less or equal to $n-2$.
\end{corollary}

\begin{proof}
Any filiform Lie algebra $\mathcal{N}$ is obtained as a deformation of the Lie
algebra $L_{n}$, since $\chi\left(  L_{n}\right)  =n-2$\, and  using Proposition \ref{propDeformation}, one has $\chi\left(  \mathcal{N}\right)  \leq n-2.$
\end{proof}

\section{Index of Filiform Lie algebras of dimension $\leq8$}

In this section, we compute the index of $n$-dimensional Filiform Lie algebras
with $n<8$. Let $\mathcal{G}$ be an $n$-dimensional Lie algebra. We set
$\left\{  x_{1},x_{2},......,x_{n}\right\}  $ be a fixed basis of
$\mathbb{V}=\mathcal{G}$, $\left\{  x_{1}^{\ast},x_{2}^{\ast},...,x_{n}^{\ast
}\right\}  $ and $f=\sum_{i\geq1}{p_{i}x_{i}^{\ast}}$.

\subsection{Filiform Lie algebras of dimension less than 6}

Any $n$-dimensional Lie algebras with $n<5$ is isomorphic to one of the
following Lie algebras.

\textbf{Dimension 1 and 2} We have only the abelian Lie algebras.

\textbf{Dimension 3}

\hspace{1cm}$\mathcal{F}_{3}^{1}$ : $\left[  x_{1},x_{2}\right]  =x_{3}.$

\textbf{Dimension 4}

\hspace{1cm}$\mathcal{F}_{4}^{1}$ : $\left[  x_{1},x_{2}\right]  =x_{3}$ $,$
$\left[  x_{1},x_{3}\right]  =x_{4}.$

\textbf{Dimension 5}

\hspace{1cm}$\mathcal{F}_{5}^{1}$ : $\left[  x_{1},x_{i}\right]  =x_{i+1},$
for $i=2,3,4.$

\hspace{1cm}$\mathcal{F}_{5}^{2}$ : $\left[  x_{1},x_{i}\right]  =x_{i+1},$
for $i=2,3,4$ and $\left[  x_{2},x_{3}\right]  =x_{5}.$

\vspace{0.5cm} The computations of the index using Proposition \ref{Dixmier}
lead to the following result.

\begin{proposition}
The index of $n$-dimensional filiform Lie algebras with $n\leq5$ are
\[
\chi\left(  \mathcal{F}_{3}^{1}\right)  =1,\quad\chi\left(  \mathcal{F}%
_{4}^{1}\right)  =2,\quad\chi\left(  \mathcal{F}_{5}^{1}\right)  =3,\quad
\chi\left(  \mathcal{F}_{5}^{2}\right)  =1.
\]
The regular vectors of $\mathcal{F}_{n}^{1}$ for $n=3,4,5$ are of the form
$f=\sum_{i=1}^{5}{p_{i}x_{i}^{\ast}}$ with one of $p_{i}\neq0$ i$\in
\{3,4,5\}$\newline The regular vectors of $\mathcal{F}_{5}^{2}$ are of the
form $f=\sum_{i=1}^{5}{p_{i}x_{i}^{\ast}}$with $p_{i}\neq0.i\in\{3,4,5\}$%
\newline
\end{proposition}

\begin{proof}
The filiform Lie algebras $\mathcal{F}_{3}^{1}$, $\mathcal{F}_{4}^{1}$ and
$\mathcal{F}_{5}^{1}$ are of type $L_{n}$. For $\mathcal{F}_{5}^{2}$, the
corresponding matrix is of rank 4, then the index is one. The regular vector
are obtained by direct calculation.
\end{proof}

\subsection{Filiform Lie algebras of dimension 6}

Any $n$-dimensional Lie algebras with $n=6$ is isomorphic to one of the
following Lie algebras.

$\mathcal{F}_{6}^{1}:$ $\left[  x_{1},x_{i}\right]  =x_{i+1},$ for $i=2,3,4,5$

$\mathcal{F}_{6}^{2}:$ $\left[  x_{1},x_{i}\right]  =x_{i+1},$ for
$i=2,3,4,5$, $\left[  x_{2},x_{3}\right]  =x_{6}$

$\mathcal{F}_{6}^{3}:$ $\left[  x_{1},x_{i}\right]  =x_{i+1},$ for
$i=2,3,4,5$, $\left[  x_{2},x_{5}\right]  =x_{6}$, and $\left[  x_{3}%
,x_{4}\right]  =-x_{6}$

$\mathcal{F}_{6}^{4}:$ $\left[  x_{1},x_{i}\right]  =x_{i+1},$ for
$i=2,3,4,5$, $\left[  x_{2},x_{3}\right]  =x_{5}$, and $\left[  x_{2}%
,x_{4}\right]  =x_{6}$

$\mathcal{F}_{6}^{5}:$ $\left[  x_{1},x_{i}\right]  =x_{i+1},$ for
$i=2,3,4,5$, $\left[  x_{2},x_{3}\right]  =x_{5}-x_{6},\left[  x_{2}%
,x_{4}\right]  =x_{6},\left[  x_{2},x_{5}\right]  =x_{6},\left[  x_{3}%
,x_{4}\right]  =-x_{6}$.
\begin{proposition}
The index of $6$-dimensional filiform Lie algebras are%

\begin{align*}
\chi\left(  \mathcal{F}_{6}^{i}\right)   &  =2\text{ for }i=2,4,3,5.\\
\chi\left(  \mathcal{F}_{6}^{1}\right)   &  =4
\end{align*}

The regular vectors of $\mathcal{F}_{6}^{1}$ are of the form $f=\sum_{i=1}%
^{6}p_{i}{x_{i}^{\ast}}$with one of $p_{i}\neq0.i=\left\{  3,...,6\right\}
$(class of $L_{n}$ algebra ). 

The regular vectors of $\mathcal{F}%
_{6}^{2}$ are of the form $f=p_{1}x_{1}^{\ast}+p_{2}x_{2}^{\ast}+p({x_{3}%
^{\ast}+x_{4}^{\ast}+x_{5}^{\ast}})+$ $p_{5}x_{5}^{\ast}$.\ 

The regular
vectors of $\mathcal{F}_{6}^{4}$ are of the form $f=\sum_{i=1}^{5}p_{i}%
{x_{i}^{\ast}}$ with $p_{6}=0.$

The regular vectors of $\mathcal{F}_{6}^{i}$ for $i=3,5$ are of the form
$f=\sum_{i=1}^{6}p_{i}{x_{i}^{\ast}}$ with one of $p_{i}\neq0$ $i\in
\{3,...6\}$.
\end{proposition}
\subsection{Filiform Lie algebras of dimension 7}

Any $n$-dimensional Lie algebras with $n=7$ is isomorphic to one of the
following Lie algebras.

$\mathcal{F}_{7}^{1}:$ $\left[  x_{1},x_{i}\right]  =x_{i+1},$ for
$i=2,3,4,5,$ \  \  $\left[  x_{1},x_{6}\right]  =\alpha x_{7},$
$\ \  \left[  x_{2},x_{3}\right]  =\left(  1+\alpha\right)  x_{5},$
$\ \ \left[  x_{2},x_{4}\right]  =\left(  1+\alpha\right)  x_{6},$ $\ \left[
x_{3},x_{4}\right]  =x_{7}.$

$\mathcal{F}_{7}^{2}:$ $\left[  x_{1},x_{i}\right]  =x_{i+1},$ for
$i=2,3,4,5,6,$ $\ \left[  x_{2},x_{3}\right]  =x_{5},$ $\ $%
\ $\ \  \left[  x_{2},x_{4}\right]  =x_{6},$
$\ \ \ \left[  x_{2},x_{5}\right]  =x_{7}.$

$\mathcal{F}_{7}^{3}:$ $\left[  x_{1},x_{i}\right]  =x_{i+1},$ for
$i=2,3,4,5,6,$ $\ \left[  x_{2},x_{3}\right]  =x_{5}+x_{6},$ $\ \left[
x_{2},x_{4}\right]  =x_{6},$ $\ \ \left[  x_{2}%
,x_{5}\right]  =x_{7}.$

$\mathcal{F}_{7}^{4}:$ $\left[  x_{1},x_{i}\right]  =x_{i+1},$ for
$i=2,3,4,5,6,$ $\ \left[  x_{2},x_{3}\right]  =x_{6},$
$\ \  \left[  x_{2},x_{4}\right]  =x_{7},$
$\ \  \left[  x_{2},x_{5}\right]  =x_{7},$
$\ \  \left[  x_{3},x_{4}\right]  =-x_{7}.$

$\mathcal{F}_{7}^{5}:$ $\left[  x_{1},x_{i}\right]  =x_{i+1},$ for
$i=2,3,4,5,6,$ $\ \left[  x_{2},x_{3}\right]  =x_{6}+x_{7},$ $\ \left[
x_{2},x_{4}\right]  =x_{7}.$

$\mathcal{F}_{7}^{6}:$ $\left[  x_{1},x_{i}\right]  =x_{i+1},$ for
$i=2,3,4,5,6,$ $\ \left[  x_{2},x_{3}\right]  =x_{6},$
$\ \  \left[  x_{2},x_{4}\right]  =x_{7}.$

$\mathcal{F}_{7}^{7}:$ $\left[  x_{1},x_{i}\right]  =x_{i+1},$ for
$i=2,3,4,5,6,$ $\ \left[  x_{2},x_{3}\right]  =x_{7}$

$\mathcal{F}_{7}^{8}:$ $\left[  x_{1},x_{i}\right]  =x_{i+1},$ for
$i=2,3,4,5,6 $ \ (class of $L_{n}$ algebra ).

\begin{proposition}
The index of $7$-dimensional filiform Lie algebras are%
\begin{align*}
\chi\left(  \mathcal{F}_{7}^{i}\right)   &  =3\text{ for }%
i=2,3,5,6,7\ \ \ \chi\left(  \mathcal{F}_{7}^{4}\right)  =1,\\
\chi\left(  \mathcal{F}_{7}^{1}\right)   &  =
\begin{cases}
1\text{ if }\alpha\neq\left\{  0,-1\right\}, \\
3\text{ if }\alpha=0.
\end{cases}
 \\
\chi\left(  \mathcal{F}_{7}^{8}\right)   &  =5.
\end{align*}

The regular vectors of $\mathcal{F}_{7}^{i}$ are given by the following table

\begin{center}%
\begin{tabular}
[c]{|c|c|}\hline
item & regular vectors\\\hline
$i=1$ & $f=\sum_{i=1}^{5}p_{i}{x_{i}^{\ast}+p(x_{6}^{\ast}+x_{7}^{\ast}%
)},\text{ if }\alpha=0$\\
\  & $f=\sum_{i=1}^{6}p_{i}{x_{i}^{\ast}}\text{ with }p_{i}\neq0,\text{ if
}\alpha\neq0$\\\hline
$i=2$ & $f=p_{1}x_{1}^{\ast}+p_{2}x_{2}^{\ast}+p(x_{4}^{\ast}+x_{5}^{\ast
}+x_{6}^{\ast})\text{ with }p\neq0$\\\hline
$i=3$ & $f=\sum_{i=1}^{4}p_{i}{x_{i}^{\ast}}$\\\hline
$i=4$ & $f=\sum_{i=1}^{7}p_{i}{x_{i}^{\ast}}\text{ {with }}p_{4}=0,p_{3}%
=0$\\\hline
$i=5$ & $f=p_{1}x_{1}^{\ast}+p_{2}x_{2}^{\ast}+p_{3}{x_{3}^{\ast}+p_{4}%
{x_{4}^{\ast}+}p(x_{5}^{\ast}+x_{6}^{\ast}})$\\\hline
$i=6$ & $f=p_{1}x_{1}^{\ast}+p_{2}x_{2}^{\ast}+p_{3}{x_{3}^{\ast}%
+p(x_{4}^{\ast}+x_{5}^{\ast}})$\\\hline
$i=7$ & $f=p_{1}x_{1}^{\ast}+p_{2}x_{2}^{\ast}+p_{3}{x_{3}^{\ast}+p_{4}%
{x_{4}^{\ast}+}p(x_{6}^{\ast}+x_{7}^{\ast}})$\\\hline
$i=8$ & $f=\sum_{i=1}^{7}p_{i}{x_{i}^{\ast}}\text{with one of }p_{i}%
\neq0\ \ i\in\{3,...7\}$\\\hline
\end{tabular}

\end{center}
\end{proposition}
\section{Index of Graded quasi-filiform Lie algebras}

The classification of naturally graded quasi-filiform Lie algebras is known and given in
\cite{Gomez2}. They have the characteristic sequence $\left(  n-2,1,1\right)
$ where $n$ is the dimension of the algebra.

\begin{definition}
\cite{Gomez2} An $n$-dimensional nilpotent Lie algebra $\mathcal{G}$ is said to
be quasi-filiform if $C^{n-3}\mathcal{G}\neq0$ and $C^{n-2}\mathcal{G}=0,$
where $C^{0}\mathcal{G}=\mathcal{G},$ $C^{i}\mathcal{G}=\left[  \mathcal{G}%
,C^{i-1}\mathcal{G}\right]  .$
\end{definition}

In the following we describe the classification of naturally quasi-graded
filiform Lie algebras

let $\mathcal{B}=\left\{  x_{0}x_{2},...,x_{n-1}\right\}  $ be a basis of
$\mathcal{G}:$

\subsection{Naturally graded Quasi-filiform Lie algebras}

We consider the following classes of $n$-dimensional Lie algebras which are
naturally graded quasi-filiform Lie algebras,

we set

Split : $L_{n-1}\oplus\mathbb{C}$\ $\left(  n\geq4\right)  :$

\hspace{1cm}$\left[  x_{0},x_{i}\right]  =x_{i+1}.$\hspace{1cm}$1\leq i\leq
n-3.$

\hspace{1cm}$Q_{n-1}\oplus\mathbb{C}$ $\left(  n\geq7,\text{\ }n\text{
odd}\right) , $

\hspace{1cm}$\left[  x_{0},x_{i}\right]  =x_{i+1}.$\hspace{1cm}$1\leq i\leq
n-3,$

\hspace{1cm}$\left[  x_{i},x_{n-2-i}\right]  =\left(  -1\right)  ^{i-1}%
x_{n-2}.$\hspace{1cm}$1\leq i\leq\frac{n-3}{2}.$

Principal : $\mathcal{L}_{\left(  n,r\right)  }$\ $\left(  n\geq5,\ r\text{
odd},\ 3\leq r\leq2\left[  \frac{n-1}{2}\right]  -1\right)  :$

\hspace{1cm}$\left[  x_{0},x_{i}\right]  =x_{i+1},$\hspace{1cm}$1\leq i\leq
n-3,$

\hspace{1cm}$\left[  x_{i},x_{r-i}\right]  =\left(  -1\right)  ^{i-1}x_{n-1},$\hspace{1cm}$1\leq i\leq\frac{r-1}{2},$

\hspace{1cm}$Q_{\left(  _{n,r}\right)  }$ $\left(  n\geq7,\ \ n\text{
odd},\ r\text{ odd},3\leq r\leq n-4\right)  $:

$\hspace{1cm}\left[  x_{0},x_{i}\right]  =x_{i+1},$\hspace{1cm}$1\leq i\leq
n-3,$

\hspace{1cm}$\left[  x_{i},x_{r-i}\right]  =\left(  -1\right)  ^{i-1}x_{n-1}%
,$\hspace{1cm}$1\leq i\leq\frac{r-1}{2},$

\hspace{1cm}$\left[  x_{i},x_{n-2-i}\right]  =\left(  -1\right)  ^{i-1}%
x_{n-2},$\hspace{1cm}$1\leq i\leq\frac{n-3}{2}.$

Terminal : $\mathcal{T}_{\left(  n,n-3\right)  }$ $\left(  n\text{ even}%
,n\geq6\right)  :$

\hspace{1cm}$\left[  x_{0},x_{i}\right]  =x_{i+1},$\hspace{1cm}$1\leq i\leq
n-3,$

\hspace{1cm}$\left[  x_{n-1},x_{1}\right]  =\frac{n-4}{2}x_{n-2},$

\hspace{1cm}$\left[  x_{i},x_{n-3-i}\right]  =\left(  -1\right)  ^{i-1}\left(
x_{n-3}+x_{n-1}\right)  ,$\hspace{1cm}$1\leq i\leq\frac{n-4}{2},$

$\hspace{1cm}\left[  x_{i},x_{n-2-i}\right]  =\left(  -1\right)  ^{i-1}%
\frac{n-2-2i}{2}x_{n-2},$\hspace{1cm}$1\leq i\leq\frac{n-4}{2},$

\hspace{1cm}$\mathcal{T}_{\left(  n,n-4\right)  }$\ $\left(  n\text{
odd},n\geq7\right) , $

\hspace{1cm}$\left[  x_{0},x_{i}\right]  =x_{i+1},$\hspace{1cm}$1\leq i\leq
n-3,$

\hspace{1cm}$\left[  x_{n-1},x_{1}\right]  =\frac{n-5}{2}x_{n-4+i},$%
\hspace{1cm}$1\leq i\leq2,$

\hspace{1cm}$\left[  x_{i},x_{n-4-i}\right]  =\left(  -1\right)  ^{i-1}\left(
x_{n-4}+x_{n-1}\right)  ,$\hspace{1cm}$1\leq i\leq\frac{n-5}{2},$

\hspace{1cm}$\left[  x_{i},x_{n-3-i}\right]  =\left(  -1\right)  ^{i-1}%
\frac{n-3-2i}{2}x_{n-2},$\hspace{1cm}$1\leq i\leq\frac{n-5}{2},$

\hspace{1cm}$\left[  x_{i},x_{n-2-i}\right]  =\left(  -1\right)  ^{i-1}\left(
i-1\right)  \frac{n-3-i}{2}x_{n-2},$\hspace{1cm}$1\leq i\leq\frac{n-3}{2}.$

Moreover, we have the following $7$-dimensional and $9$-dimensional Lie
algebra \cite{GarciaVergnolle}.

$\varepsilon_{\left(  7,3\right)  }:
\begin{cases}
\left[  x_{0},x_{i}\right]  =x_{i+1},\ \ 1\leq i\leq4,\\
\left[  x_{6},x_{i}\right]  =x_{3+i},\ \ 1\leq i\leq2,\\
\left[  x_{1},x_{2}\right]  =x_{3}+x_{6},\ \ \\
\left[  x_{1},x_{i}\right]  =x_{i+1},\ \ 3\leq i\leq4.
\end{cases}
 $

$\varepsilon_{\left(  9,5\right)  }^{1}:
\begin{cases}
\left[  x_{0},x_{i}\right]  =x_{i+1},\ \ 1\leq i\leq6,\\
\left[  x_{8},x_{i}\right]  =2x_{5+i},\ \ 1\leq i\leq2,\\
\left[  x_{1},x_{4}\right]  =x_{5}+x_{8},\ \ \left[  x_{1}%
,x_{5}\right]  =2x_{6},\\
\left[  x_{1},x_{6}\right]  =3x_{7},\ \ \left[  x_{2},x_{3}\right]
=-x_{5}-x_{8},\\
\left[  x_{2},x_{4}\right]  =-x_{6},\ \ \left[  x_{2},x_{5}\right]
=-x_{7}.\hspace{1cm}%
\end{cases}  $

$\varepsilon_{\left(  9,5\right)  }^{2}:
\begin{cases}\left[  x_{0},x_{i}\right]  =x_{i+1},\ \ 1\leq i\leq6,\\
\left[  x_{8},x_{i}\right]  =2x_{5+i},\ \ 1\leq i\leq2,\\
\left[  x_{1},x_{4}\right]  =x_{5}+x_{8},\ \ \left[
x_{1},x_{5}\right]  =2x_{6},\\
\left[  x_{1},x_{6}\right]  =x_{7},\ \ \left[  x_{2}%
,x_{3}\right]  =-x_{5}-x_{8},\\
\left[  x_{2},x_{4}\right]  =-x_{6},\left[  x_{2},x_{5}\right]  =x_{7},\left[
x_{3},x_{4}\right]  =-2x_{7}.
\end{cases}  $

$\varepsilon_{\left(  9,5\right)  }^{3}:
\begin{cases}\left[  x_{0},x_{i}\right]  =x_{i+1},\ \ 1\leq i\leq6,\\
\left[  x_{0},x_{8}\right]  =x_{6},\ \ 1\leq i\leq2,\\
\left[  x_{1},x_{4}\right]  =x_{8},\ \ \left[  x_{3}%
,x_{4}\right]  =-3x_{7},\\
\left[  x_{2},x_{4}\right]  =-x_{6},\ \ \left[  x_{1}%
,x_{5}\right]  =2x_{6},\\
\left[  x_{2},x_{3}\right]  =-x_{8},\ \ \left[  x_{2}%
,x_{5}\right]  =2x_{7}.
\end{cases} $

are graded quasi-filiform Lie algebras.

\begin{theorem}
\cite{Gomez2} Every $n$-dimensional naturally graded quasi-filiform Lie
algebra is isomorphic to one of the following Lie algebras :

\begin{itemize}
\item If $n$ is even to $L_{n-1}\oplus\mathbb{C}$, $\mathcal{T}_{\left(
n,n-3\right)  }$, or $\mathcal{L}_{\left(  n,r\right)  }$ with $r$ odd and
$3\leq r\leq n-3.$

\item If $n$ is odd to $L_{n-1}\oplus\mathbb{C}$, $Q_{n-1}\oplus\mathbb{C}$,
$\mathcal{L}_{\left(  n,n-2\right)  }$, $\mathcal{T}_{\left(  n,n-4\right)  }%
$, $\mathcal{L}_{\left(  n,r\right)  }$, or $Q_{\left(  _{n,r}\right)  }$
with\ $r$ odd, and\ $3\leq r\leq n-4$. In the case of $n=7$ and $n=9,$ we add
$\varepsilon_{\left(  7,3\right)  },$ $\varepsilon_{\left(  9,5\right)  }%
^{1},$ $\varepsilon_{\left(  9,5\right)  }^{2},$ $\varepsilon_{\left(
9,5\right)  }^{3}$.
\end{itemize}
\end{theorem}

\subsubsection{Index of graded quasi-filiform Lie algebras\ :}

In the following we compute the index of graded quasi--filiform Lie algebras.
Let $\mathcal{G}$ be a $n$-dimensional graded quasi-filiform Lie algebra

\begin{theorem}
Index of graded quasi-filiform Lie algebras are

\textbf{case where }$n$\textbf{\ is even }

\begin{enumerate}
\item $\chi(L_{n-1}\oplus\mathbb{C})=n-2.$

\item $\chi(\mathcal{T}_{\left(  n,n-3\right)  })=2.$

\item $\chi(\mathcal{L}_{\left(  n,r\right)  })=n-r-1, \ \ 3\leq
r\leq n-3.$
\end{enumerate}

\textbf{case where }$n$\textbf{\ is odd\ :}

\begin{enumerate}
\item $\chi(L_{n-1}\oplus\mathbb{C})=n-2.$

\item $\chi(Q_{n-1}\oplus\mathbb{C})=3.$

\item $\chi(\mathcal{L}_{\left(  n,n-2\right)  })=3.$

\item $\chi(\mathcal{T}_{\left(  n,n-4\right)  })=3.$

\item $\chi(\mathcal{L}_{\left(  n,r \right)  })=n-r-1, \ Ê\ 3 \leq r \leq n-3.$

\item $\chi(Q_{\left(  n,r\right)  })=3.$

\item $\chi\left(  \varepsilon_{\left(  7,3\right)  }\right)  =3.$

\item $\chi\left(  \varepsilon_{\left(  9,5\right)  }^{1}\right)  =3.$

\item $\chi\left(  \varepsilon_{\left(  9,5\right)  }^{i}\right)
=2, \ \ \  i=2,3.$
\end{enumerate}
\end{theorem}

\begin{proof}
\textbf{case where }$n$\textbf{\ is even }

The corresponding matrix to the graded quasi-filiform Lie algebra
$L_{n-1}\oplus\mathbb{C}$ is of the form

$$\left(
\begin{tabular}
[c]{llllll}%
$0$ & $x_{2}$ & $\cdots$ & $x_{n-1}$ & $0$ & $0$\\
$-x_{2}$ & $0$ & $\cdots$ & $0$ & $0$ & $0$\\
$\vdots$ & $\vdots$ & $\vdots$ & $\vdots$ & $\vdots$ & $\vdots$\\
$-x_{n-1}$ & $0$ & $\cdots$ & $0$ & $0$ & $0$\\
$0$ & $0$ & $\cdots$ & $0$ & $0$ & $0$\\
$0$ & $0$ & $\cdots$ & $0$ & $0$ & $0$%
\end{tabular}
\right)  $$\\

Its rank is 2, then $\chi(L_{n-1}\oplus\mathbb{C})=n-2.$

The corresponding matrix to the graded quasi-filiform Lie algebra
$\mathcal{T}_{\left(  n,n-3\right)  }$ is of the form

$$\left(
\begin{tabular}
[c]{llllllll}%
$0$ & $x_{2}$ & $x_{3}$ & $\cdots$ & $x_{n-3}$ & $x_{n-2}$ & $0$ & $0$\\
$-x_{2}$ & $0$ & $0$ & $\cdots$ & $x_{n-3}+x_{n-1}$ & $\left(  \frac{n-4}%
{2}\right)  x_{n-2}$ & $0$ & $\left(  \frac{n-4}{2}\right)  x_{n-2}$\\
$-x_{3}$ & $0$ & $0$ & $\cdots$ & $-\left(  \frac{n-6}{2}\right)  x_{n-2}$ &
$0$ & $0$ & $0$\\
$\vdots$ & $\vdots$ & $\vdots$ & $\vdots$ & $\vdots$ & $\vdots$ & $\vdots$ &
$\vdots$\\
$-x_{n-3}$ & $-x_{n-3}-x_{n-1}$ & $\left(  \frac{n-6}{2}\right)  x_{n-2}$ &
$\cdots$ & $0$ & $0$ & $0$ & $0$\\
$-x_{n-2}$ & $-\left(  \frac{n-4}{2}\right)  x_{n-2}$ & $0$ & $\cdots$ & $0$ &
$0$ & $0$ & $0$\\
$0$ & $0$ & $0$ & $\cdots$ & $0$ & $0$ & $0$ & $0$\\
$0$ & $-\left(  \frac{n-4}{2}\right)  x_{n-2}$ & $0$ & $\cdots$ & $0$ & $0$ &
$0$ & $0$%
\end{tabular}
\right)  $$\\

Its rank is $n-2$, then $\chi(\mathcal{T}_{\left(  n,n-3\right)  })=2.$

The corresponding matrix to the graded quasi-filiform Lie algebra
$\mathcal{L}_{\left(  n,r\right)  }$ is of the form

$$\left(
\begin{array}
[c]{cccccccccc}%
0 & x_{2} & x_{3} & \cdots & x_{r} & \cdots & x_{n-3} & x_{n-2} & 0 & 0\\
-x_{2} & 0 & 0 & \cdots & -x_{n-1} & \cdots & 0 & 0 & 0 & 0\\
-x_{3} & 0 & 0 & \cdots & 0 & \cdots & 0 & 0 & 0 & 0\\
\vdots & \vdots & \vdots & \vdots & \vdots & \vdots & \vdots & \vdots & \vdots
& \vdots\\
-x_{r} & x_{n-1} & 0 & \cdots & 0 & \cdots & 0 & 0 & 0 & 0\\
\vdots & \vdots & \vdots & \vdots & \vdots & \cdots & \vdots & \vdots & \vdots
& \vdots\\
-x_{n-3} & 0 & 0 & \cdots & 0 & \cdots & 0 & 0 & 0 & 0\\
-x_{n-2} & 0 & 0 & \cdots & 0 & \cdots & 0 & 0 & 0 & 0\\
0 & 0 & 0 & \cdots & 0 & \cdots & 0 & 0 & 0 & 0\\
0 & 0 & 0 & \cdots & 0 & \cdots & 0 & 0 & 0 & 0
\end{array}
\right)  $$\\

For $3\leq r\leq n-3$, its rank is $r+1$. Then $\chi(\mathcal{L}_{\left(
n,r\right)  })=n-r-1.$

\textbf{case where }$n$\textbf{\ is odd :}

The corresponding matrix to the graded quasi-filiform Lie algebra
$L_{n-1}\oplus\mathbb{C}$ is of the form

$$\left(
\begin{tabular}
[c]{llllll}%
$0$ & $x_{2}$ & $\cdots$ & $x_{n-1}$ & $0$ & $0$\\
$-x_{2}$ & $0$ & $\cdots$ & $0$ & $0$ & $0$\\
$\vdots$ & $\vdots$ & $\vdots$ & $\vdots$ & $\vdots$ & $\vdots$\\
$-x_{n-1}$ & $0$ & $\cdots$ & $0$ & $0$ & $0$\\
$0$ & $0$ & $\cdots$ & $0$ & $0$ & $0$\\
$0$ & $0$ & $\cdots$ & $0$ & $0$ & $0$%
\end{tabular}
\right)  $$\\

Its rank is 2, then $\chi(L_{n-1}\oplus\mathbb{C})=n-2.$

The corresponding matrix to the graded quasi-filiform Lie algebra
$Q_{n-1}\oplus\mathbb{C}$ is of the form

$$\left(
\begin{tabular}
[c]{llllllll}%
$0$ & $x_{2}$ & $x_{3}$ & $\cdots$ & $x_{n-3}$ & $x_{n-2}$ & $0$ & $0$\\
$-x_{2}$ & $0$ & $0$ & $\cdots$ & $0$ & $-x_{n-2}$ & $0$ & $0$\\
$-x_{3}$ & $0$ & $0$ & $\cdots$ & $-x_{n-2}$ & $0$ & $0$ & $0$\\
$\vdots$ & $\vdots$ & $\vdots$ & $\vdots$ & $\vdots$ & $\vdots$ & $\vdots$ &
$\vdots$\\
$-x_{n-3}$ & $0$ & $_{-}x_{n-2}$ & $\cdots$ & $0$ & $0$ & $0$ & $0$\\
$-x_{n-2}$ & $x_{n-2}$ & $0$ & $\cdots$ & $0$ & $0$ & $0$ & $0$\\
$0$ & $0$ & $0$ & $\cdots$ & $0$ & $0$ & $0$ & $0$\\
$0$ & $0$ & $0$ & $\cdots$ & $0$ & $0$ & $0$ & $0$%
\end{tabular}
\right)  $$\\

Its rank is $n-3$, then $\chi(Q_{n-1}\oplus\mathbb{C})=3.$

The corresponding matrix to the graded quasi-filiform Lie algebra
$\mathcal{L}_{\left(  n,n-2\right)  }$ is of the form

$$\left(
\begin{tabular}
[c]{lllllllll}%
$0$ & $x_{2}$ & $x_{3}$ & $x_{4}$ & $\cdots$ & $x_{n-3}$ & $x_{n-2}$ & $0$ &
$0$\\
$-x_{2}$ & $0$ & $0$ & $0$ & $\cdots$ & $0$ & $-x_{n-1}$ & $0$ & $0$\\
$-x_{3}$ & $0$ & $0$ & $0$ & $\cdots$ & $x_{n-1}$ & $0$ & $0$ & $0$\\
$-x_{4}$ & $0$ & $0$ & $0$ & $\cdots$ & $0$ & $0$ & $0$ & $0$\\
$\vdots$ & $\vdots$ & $\vdots$ & $\vdots$ & $\vdots$ & $0$ & $0$ & $0$ & $0$\\
$-x_{n-3}$ & $0$ & $-x_{n-1}$ & $0$ & $\cdots$ & $0$ & $0$ & $0$ & $0$\\
$-x_{n-2}$ & $x_{n-1}$ & $0$ & $0$ & $\cdots$ & $0$ & $0$ & $0$ & $0$\\
$0$ & $0$ & $0$ & $0$ & $\cdots$ & $0$ & $0$ & $0$ & $0$\\
$0$ & $0$ & $0$ & $0$ & $\cdots$ & $0$ & $0$ & $0$ & $0$%
\end{tabular}
\right)  $$\\

Its rank is $n-3$, then $\chi(\mathcal{L}_{\left(  n,n-2\right)  })=3.$

The corresponding matrix to the graded quasi-filiform Lie algebra
$\mathcal{L}_{\left(  n,r\right)  }$ is of the form

$$\left(
\begin{array}
[c]{cccccccccc}%
0 & x_{2} & x_{3} & \cdots & x_{r} & \cdots & x_{n-3} & x_{n-2} & 0 & 0\\
-x_{2} & 0 & 0 & \cdots & -x_{n-1} & \cdots & 0 & 0 & 0 & 0\\
-x_{3} & 0 & 0 & \cdots & 0 & \cdots & 0 & 0 & 0 & 0\\
\vdots & \vdots & \vdots & \vdots & \vdots & \vdots & \vdots & \vdots & \vdots
& \vdots\\
-x_{r} & x_{n-1} & 0 & \cdots & 0 & \cdots & 0 & 0 & 0 & 0\\
\vdots & \vdots & \vdots & \vdots & \vdots & \cdots & \vdots & \vdots & \vdots
& \vdots\\
-x_{n-3} & 0 & 0 & \cdots & 0 & \cdots & 0 & 0 & 0 & 0\\
-x_{n-2} & 0 & 0 & \cdots & 0 & \cdots & 0 & 0 & 0 & 0\\
0 & 0 & 0 & \cdots & 0 & \cdots & 0 & 0 & 0 & 0\\
0 & 0 & 0 & \cdots & 0 & \cdots & 0 & 0 & 0 & 0
\end{array}
\right)  $$\\

Its rank is $r+1$, then $\chi(\mathcal{L}_{\left(  n,r\right)  })=n-r-1$,
$3\leq r\leq n-4$.

The corresponding matrix to the graded quasi-filiform Lie algebra
 $\mathcal{T}_{\left(  n,n-4\right)  }$ is of the form\\

\begin{small}
$$\left(
\begin{array}
[c]{llllllllll}%
0 & x_{2} & x_{3} & x_{4} & \cdots & x_{n-4} & x_{n-3} & x_{n-2} & 0 & 0\\
-x_{2} & 0 & 0 & 0 & \cdots & x_{n-4}+x_{n-1} & \frac{n-5}{2}x_{n-3} & 0 & 0 &
-\left(  \frac{n-5}{2}\right)  x_{n-3}\\
-x_{3} & 0 & 0 & 0 & \cdots & -\left(  \frac{n-7}{2}\right)  x_{n-3} &
\frac{n-5}{2}x_{n-2} & 0 & 0 & -\left(  \frac{n-5}{2}\right)  x_{n-2}\\
-x_{4} & 0 & 0 & 0 & \vdots & -2\left(  \frac{n-6}{2}\right)  x_{n-2} & 0 &
0 & 0 & 0\\
\vdots & \vdots & \vdots & \vdots & \cdots & \vdots & \vdots & \vdots & \vdots
& \vdots\\
-x_{n-4} & -x_{n-4}-x_{n-1} & \left(  \frac{n-7}{2}\right)  x_{n-3} & 2\left(
\frac{n-6}{2}\right)  x_{n-2} & \cdots & 0 & 0 & 0 & 0 & 0\\
-x_{n-3} & -\left(  \frac{n-5}{2}\right)  x_{n-3} & -\left(  \frac{n-5}%
{2}\right)  x_{n-2} & 0 & \cdots & 0 & 0 & 0 & 0 & 0\\
-x_{n-2} & 0 & 0 & 0 & \cdots & 0 & 0 & 0 & 0 & 0\\
0 & 0 & 0 & 0 & \cdots & 0 & 0 & 0 & 0 & 0\\
0 & \left(  \frac{n-5}{2}\right)  x_{n-3} & \left(  \frac{n-5}{2}\right)
x_{n-2} & 0 & \cdots & 0 & 0 & 0 & 0 & 0
\end{array}
\right)  $$
\end{small}\\

Its rank \ is $n-3$,  then $\chi\left(  \mathcal{T}_{\left(  n,n-4\right)
}\right)  =3.$

The corresponding matrix to the graded quasi-filiform Lie algebra $Q_{\left(
_{n,r}\right)  }$ is of the form

$$\left(
\begin{tabular}
[c]{llllllllll}%
$0$ & $x_{2}$ & $x_{3}$ & $\cdots$ & $x_{r}$ & $\cdots$ & $x_{n-3}$ &
$x_{n-2}$ & $0$ & $0$\\
-$x_{2}$ & $0$ & $0$ & $\cdots$ & $-x_{n-1}$ & $\cdots$ & $0$ & $-x_{n-2}$ &
$0$ & $0$\\
$-x_{3}$ & $0$ & $0$ & $\cdots$ & $0$ & $\cdots$ & $0$ & $0$ & $0$ & $0$\\
$\vdots$ & $\vdots$ & $\vdots$ & $\cdots$ & $\vdots$ & $\vdots$ & $\vdots$ &
$\vdots$ & $\vdots$ & $\vdots$\\
$-x_{r}$ & $x_{n-1}$ & $0$ & $\cdots$ & $0$ & $\cdots$ & $0$ & $0$ & $0$ &
$0$\\
$\vdots$ & $\vdots$ & $\vdots$ & $\cdots$ & $\vdots$ & $\vdots$ & $\vdots$ &
$\vdots$ & $\vdots$ & $\vdots$\\
$-x_{n-3}$ & $0$ & $0$ & $\cdots$ & $0$ & $\cdots$ & $0$ & $0$ & $0$ & $0$\\
$-x_{n-2}$ & $x_{n-2}$ & $0$ & $\cdots$ & $0$ & $\cdots$ & $0$ & $0$ & $0$ &
$0$\\
$0$ & $0$ & $0$ & $\cdots$ & $0$ & $\cdots$ & $0$ & $0$ & $0$ & $0$\\
$0$ & $0$ & $0$ & $\cdots$ & $0$ & $\cdots$ & $0$ & $0$ & $0$ & $0$%
\end{tabular}
\right)  $$\\

For $3\leq r\leq n-4$
its rank is $n-3$.  Then $\chi(Q_{\left(  n,r\right)  })=3.$
\end{proof}

\begin{remark}
There are no Frobenius quasi-filiform Lie algebra.
\end{remark}

\subsubsection{Regular vectors}

\begin{proposition}
The  regular vectors of the families $\mathcal{T}_{\left(  n,n-3\right)  }$, $\mathcal{T}_{\left(  n,n-4\right)  }$, $\mathcal{L}_{\left(  n,r\right)  }$ and  $Q_{\left(  _{n,r}\right)  }$  are given by
the following functionals $f$ where $x_{i}^{\ast}$ are the element of the dual basis and $p_i$ are parameters. 
\begin{enumerate}
\item $\mathcal{T}_{\left(  n,n-3\right)  }:$%
\[
f=\sum_{i=0}^{n-1}{p_{i\text{ }}x_{i}^{\ast}}\text{ with }p_{n-2}\neq0.
\]

\item $\mathcal{T}_{\left(  n,n-4\right)  }:$%
\[
f=\sum_{i=0}^{n-1}{p_{i}x_{i}^{\ast}}\text{ with }p_{n-2}\neq0.
\]

\item $\mathcal{L}_{\left(  n,r\right)  }:$ $n$ odd or even and $r<n-2$ :%
\[
f=\sum_{i=0}^{n-1}{p_{i}x_{i}^{\ast}}\text{ with }p_{n-1}\neq0\text{ and one
of }{p_{i}}\neq0\text{ where }i\in\left\{  {r+1},...,{n-2}\right\}.
\]

\item $Q_{\left(  _{n,r}\right)  }:$%
\[
f=\sum_{i=0}^{n-1}{p_{i\text{ }}x_{i}^{\ast}}\text{ with }p_{n-2}\neq0.
\]

\item $\mathcal{L}_{\left(  n,n-2\right)  }:$%
\[
f=\sum_{i=0}^{n-1}{p_{i}x_{i}^{\ast}}\text{ with }p_{n-1}\neq0.
\]

\end{enumerate}
\end{proposition}

\begin{proof}
$\mathcal{T}_{\left(  n,n-3\right)  }:$

The associate system of the graded quasi-filiform Lie algebra
$\mathcal{T}_{\left(  n,n-3\right)  }$ is of the form:

$
\begin{cases}
\overset{n-3}{\underset{i=1}{\sum}}a_{i}p_{i+1}=0,\\
a_{0}p_{2}-a_{n-4}(p_{n-3}+p_{n-1})-\frac{n-4}{2}a_{n-3}p_{n-2}+\frac{n-4}%
{2}a_{n-1}p_{n-2}=0,\\
a_{0}p_{i+1}+\left(  -1\right)  ^{i}a_{n-3-i}(p_{n-3}+p_{n-1})-\left(
-1\right)  ^{i}\frac{n-2-2i}{2}a_{n-2-i}p_{n-2}=0,\ \ i=2,...n-4,\\
a_{0}p_{n-2}+\frac{n-4}{2}a_{1}p_{n-2}=0,\\
-\frac{n-4}{2}a_{1}p_{n-2}=0.
\end{cases}  $

It turns out that $p_{n-2}\neq 0$ gives a solution of this system such that $\dim\mathcal{G}^{f}%
=\chi_{\mathcal{G}}$,
then the regular vectors are given by : $f=\sum_{i=0}^{n-1}{p_{i\text{ }}%
x_{i}^{\ast}}$ with $p_{n-2}\neq0.$

$\mathcal{T}_{\left(  n,n-4\right)  }$:

The associate system is of the form:

$
\begin{cases}
\overset{n-3}{\underset{i=1}{\sum}}a_{i}p_{i+1}=0,\\
a_{0}p_{2}-a_{n-5}(p_{n-4}+p_{n-1})-\frac{n-5}{2}a_{n-4}p_{n-3}+\frac{n-5}
{2}a_{n-1}p_{n-3}=0,\\
a_{0}p_{3}+a_{n-6}(p_{n-4}+p_{n-1})+\frac{n-7}{2}a_{n-5}p_{n-3}-\frac{n-5}
{2}a_{n-4}p_{n-2}+\frac{n-5}{2}a_{n-1}p_{n-2}=0,\\
a_{0}p_{i+1}+\left(  -1\right)  ^{i}a_{n-4-i}(p_{n-4}+p_{n-1})+\left(
-1\right)  ^{i}\frac{n-3-2i}{2}a_{n-3-i}p_{n-3}-\left(  -1\right)  ^{i}%
\frac{n-3-i}{2}a_{n-2-i}p_{n-2}=0\\
\hspace{3cm}i=3,...n-5,\\
a_{0}p_{n-3}=0,\\
-\frac{n-5}{2}a_{1}p_{n-3}+\frac{n-5}{2}a_{2}p_{n-2}=0,\\
a_{0}p_{n-3}+\frac{n-5}{2}a_{1}p_{n-3}+\frac{n-5}{2}a_{2}p_{n-2}=0.
\end{cases}  $

It follows that $p_{n-2}\neq 0$ gives a  solution of this system such that $\dim\mathcal{G}^{f}%
=\chi_{\mathcal{G}},$
then the regular vectors are given by : $f=\sum_{i=0}^{n-1}{p_{i\text{ }}%
x_{i}^{\ast}}$ with $p_{n-2}\neq0.$

$\mathcal{L}_{\left(  n,r\right)  }$ $n$ odd or even and $r<n-2$.

We cancel the columns $\left(  r+1\right)  $  until  $\left(  n-1\right)  $
 and the corresponding lines. We obtain the following minor

$$\left(
\begin{array}{ccccc}
0 & 0 & \cdots & 0 & -p_{n-1}\\
0 & 0 & \cdots & p_{n-1} & 0\\
\vdots & \vdots & \vdots & \vdots & \vdots\\
0 & -p_{n-1} & \cdots & 0 & 0\\
p_{n-1} & 0 & \cdots & 0 & 0
\end{array} \right) $$

It is of non-zero determinant and this leads to 
 $f=\sum_{i=0}^{n}{\ p_{i}x_{i}^{\ast},}$ with $p_{n-1}\neq0$ and one of the ${p_{i}}$ satisfies
${p_{i}}\neq0$ where $i\in\left\{ {r+1},...,{n-2}\right\}  .$

The same reasoning  and calculations are used for $Q_{\left(  _{n,r}\right)  }$ and $\mathcal{L}_{\left(
n,n-2\right)  }.$
\end{proof}

\begin{remark}
Since $L_{n-1}\oplus\mathbb{C}$ and $Q_{n-1}\oplus\mathbb{C}$ are the central
extension of $L_{n}$ and $Q_{n},$ then the regular vectors could be given using
Proposition \ref{CentralExtension}.
\end{remark}

\section{index of Lie algebras whose nilradical is $L_{n}$ or $Q_{2n}$}

Snobel and Winternitz determined the Lie algebras whose nilradical is
isomorphic to the filiform Lie algebra $L_{n}.$
In their work this algebra is denoted by $\mathfrak{n}_{n,1}$ and it  is defined with respect to 
the basis$\left\{  x_{1},...,x_{n}\right\}  $ by 
$$
\left[  x_{i},x_{n}\right]  =x_{i-1},\ \ i=2,...,n-1.
$$

\begin{theorem}
\cite{GarciaVergnolle} Let $\tau$ be a solvable Lie algebra over a field
$\mathbb{K=R}$ or $\mathbb{C}$ and having as nilradical $\mathfrak{n}_{n,1}%
$. Then it is isomorphic to one of the following Lie algebras.

\begin{enumerate}
\item If dim $\tau=n+1$, 
set $\mathcal{B}=\left\{  x_{_{1}},...,x_{_{n}},f\right\}  $ be a basis of
$\tau$.

$\tau_{n+1,1}$  defined as 

$\hspace{1cm}\left[  f,x_{i}\right]  =\left(  n-2+\beta\right)  x_{i}, \ \ i=1,...,n-1,$

$\hspace{1cm}\left[  f,x_{n}\right]  =x_{n}.$

$\tau_{n+1,2}$ defined as 

$\hspace{1cm}\left[  f,x_{i}\right]  =x_{i},  \ \ i=1,...,n-1.$

$\tau_{n+1,3}$ defined as 

$\hspace{1cm}\left[  f,x_{i}\right]  =\left(  n-i\right)  x_{i}, \ \ i=1,...,n-1,$

$\hspace{1cm}\left[  f,x_{n}\right]  =x_{n}+x_{n-1}.$

\item If dim $\tau=n+2$, 
set $\mathcal{B}=\left\{  x_{1},...,x_{n},f_{1},f_{2}\right\}  $ be a basis of
$\tau$.

$\tau_{n+2,1}$  defined as 

$\hspace{1cm}$$ \left[  f_{1},x_{i} \right] =\left( n-1-i \right)  x_{i},   i=1,...,n-1,$

$\hspace{1cm}$$ \left[  f_{2},x_{i}\right]  =x_{i},  \   i=1,...,n-1,$

$\hspace{1cm}\left[  f_{1},x_{n}\right]  =x_{n},  \   i=1,...,n-1.$
\end{enumerate}
\end{theorem}

\subsection{Index of Lie algebras whose nilradical is $\mathfrak{n}%
_{n,1}\left(  L_{n}\right)  $}

\begin{proposition}
Index of Lie algebras whose nilradical is $\mathfrak{n}_{n,1}$ are%
\begin{align*}
& \text{If dim }\tau   =n+1, \text{ then }   \chi({\tau_{n+1,i}}) =n-1,  i=1,2,3. \\
& \text{If dim }\tau   =n+2, \text{ then } \chi({\tau_{n+2,1}}) =n-2.
\end{align*}

\end{proposition}

\begin{proof}
Set dim $\tau=n+1$. The corresponding matrix to the algebra $\tau_{n+1,1}$ is of the form:

$$\left(
\begin{tabular}
[c]{llllll}%
$0$ & $0$ & $\cdots$ & $0$ & $0$ & $-\left(  n-2+\beta\right)  x_{1}$\\
$0$ & $0$ & $\cdots$ & $0$ & $0$ & $-\left(  n-2+\beta\right)  x_{2}$\\
$\vdots$ & $\vdots$ & $\vdots$ & $\vdots$ & $\vdots$ & $\vdots$\\
$0$ & $0$ & $\cdots$ & $\cdots$ & $\cdots$ & $-\left(  n-2+\beta\right)
x_{n-1}$\\
$0$ & $0$ & $\cdots$ & $0$ & $0$ & $-x_{n}$\\
$\left(  n-2+\beta\right)  x_{1}$ & $\left(  n-2+\beta\right)  x_{2}$ &
$\cdots$ & $\left(  n-2+\beta\right)  x_{n-1}$ & $x_{n}$ & $0$%
\end{tabular}
\ \ \right)  $$

Its rank is $2$, then $\chi(\tau_{n+1,1})=n-1$.

The corresponding matrix of the Lie algebra $\tau_{n+1,2}$ is of the form:

$$\left(
\begin{tabular}
[c]{llllll}%
$0$ & $0$ & $\cdots$ & $0$ & $0$ & $-x_{1}$\\
$0$ & $0$ & $\cdots$ & $0$ & $0$ & $-x_{2}$\\
$\vdots$ & $\vdots$ & $\vdots$ & $\vdots$ & $\vdots$ & $\vdots$\\
$0$ & $0$ & $\cdots$ & $0$ & $0$ & $-x_{n-1}$\\
$0$ & $0$ & $\cdots$ & $0$ & $0$ & $0$\\
$x_{1}$ & $x_{2}$ & $\cdots$ & $x_{n-1}$ & $0$ & $0$%
\end{tabular}
\ \ \right)  $$

Its rank is $2$, then $\chi(\tau_{n+1,2})=n-1$.

The corresponding matrix to the Lie algebra $\tau_{n+1,3}$ is of the form:

$$\left(
\begin{tabular}
[c]{lllllll}%
$0$ & $0$ & $\cdots$ & $0$ & $0$ & $0$ & $-\left(  n-1\right)  x_{1}$\\
$0$ & $0$ & $\cdots$ & $0$ & $0$ & $0$ & $-\left(  n-2\right)  x_{2}$\\
$\vdots$ & $\vdots$ & $\vdots$ & $\vdots$ & $\vdots$ & $\vdots$ & $\vdots$\\
$0$ & $0$ & $\cdots$ & $0$ & $0$ & $0$ & $-\left(  n-\left(  n-2\right)
\right)  x_{n-2}$\\
$0$ & $0$ & $\cdots$ & $0$ & $0$ & $0$ & $-x_{n-1}$\\
$0$ & $0$ & $\cdots$ & $0$ & $0$ & $0$ & $-x_{n}-x_{n-1}$\\
$\left(  n-1\right)  e_{1}$ & $\left(  n-2\right)  x_{2}$ & $\cdots$ &
$\left(  n-\left(  n-2\right)  \right)  x_{n-2}$ & $x_{n-1}$ & $x_{n}+x_{n-1}$
& $0$%
\end{tabular}
\ \ \right)  $$

Its rank is $2$, then $\chi(\tau_{n+1,3})=n-1$.

If dim $\tau=n+2$,
the corresponding matrix to the Lie algebra $\tau_{n+2,1}$ is of the form:

$$\left(
\begin{tabular}
[c]{lllllll}%
$0$ & $0$ & $\cdots$ & $0$ & $0$ & $\left(  n-2\right)  x_{1}$ & $x_{1}$\\
$0$ & $0$ & $\cdots$ & $0$ & $0$ & $\left(  n-3\right)  x_{2}$ & $x_{2}$\\
$\vdots$ & $\vdots$ & $\cdots$ & $\vdots$ & $\vdots$ & $\vdots$ & $\vdots$\\
$0$ & $0$ & $\cdots$ & $0$ & $0$ & $\left(  n-n\right)  x_{n-1}$ & $x_{n-1}$\\
$0$ & $0$ & $\cdots$ & $0$ & $0$ & $x_{n}$ & $0$\\
$-\left(  n-2\right)  x_{1}$ & $-\left(  n-3\right)  x_{2}$ & $\cdots$ &
$-\left(  n-n\right)  x_{n-1}$ & $-x_{n}$ & $0$ & $0$\\
$-x_{1}$ & $-x_{2}$ & $\cdots$ & $-x_{n-1}$ & $0$ & $0$ & $0$%
\end{tabular}
\ \right)  $$

Its rank is $4$, then $\chi(\tau_{n+2,1})=n-2$.
\end{proof}

\subsubsection{Regular vectors}

If dim $\tau=n+1$

\begin{enumerate}
\item $\tau_{n+1,1}$%
\[
f=\sum_{i=1}^{n}{p_{i\text{ }}x_{i}^{\ast}}\text{ with }p_{1},...,p_{n}\neq0.
\]

\item $\tau_{n+1,2}$%
\[
f=\sum_{i=1}^{n}{p_{i\text{ }}x_{i}^{\ast}}\text{ with }p_{1},...,p_{n-1}%
\neq0.
\]

\item $\tau_{n+1,3}$%
\[
f=\sum_{i=1}^{n}{p_{i\text{ }}x_{i}^{\ast}}\text{ with }p_{1},...,p_{n}\neq0.
\]

If dim $\tau=n+2$

\item $\tau_{n+2,1}$%
\[
f=\sum_{i=1}^{n}{p_{i\text{ }}x_{i}^{\ast}}\text{ with }p_{1},...,p_{n}\neq0.
\]

\end{enumerate}

\begin{proof}
Straightforward calculations following  Remark \ref{regularvector}.
\end{proof}

\subsection{Lie algebras whose nilradical is $Q_{2n}$}

\begin{proposition}
\cite{GarciaVergnolle} Any real solvable Lie algebra of dimension $2n+1$ whose
nilradical $Q_{2n}$ is isomorphic to one of the following Lie algebras :

let $\mathcal{B}=\left\{  x_{1},...,x_{2n},y\right\}  $ be a basis of $\tau$

\begin{enumerate}
\item $\tau_{2n+1}\left(  \lambda_{2}\right)  $

$\hspace{1cm}\left[  x_{1},x_{k}\right]  =x_{k+1}, \ \  2\leq k\leq2n-2,$

$\hspace{1cm}\left[  x_{k},x_{2n+1-k}\right]  =\left(  -1\right)  ^{k}x_{2n},
\ \ 2\leq k\leq n,$

$\hspace{1cm}\left[  y,x_{1}\right]  =x_{1},$

$\hspace{1cm}\left[  y,x_{k}\right]  =\left(  k-2+\lambda_{2}\right)  x_{k},
\ \ 2\leq k\leq2n-2,$

$\hspace{1cm}\left[  y,x_{2n}\right]  =\left(  2n-3+2\lambda_{2}\right)
x_{2n}.$

\item $\tau_{2n+1}\left(  2-n,\varepsilon\right)  $

$\hspace{1cm}\left[  x_{1},x_{k}\right]  =x_{k+1}, \ \ 2\leq k\leq2n-2,$

$\hspace{1cm}\left[  x_{k},x_{2n+1-k}\right]  =\left(  -1\right)  ^{k}x_{2n}
\ \ 2\leq k\leq n,$

$\hspace{1cm}\left[  y,x_{1}\right]  =x_{1}+\varepsilon x_{2n},\  \ \varepsilon=-1,0,1,$

$\hspace{1cm}\left[  y,x_{k}\right]  =\left(  k-n\right)  x_{k}, \ \ 2\leq k\leq2n-1,$

$\hspace{1cm}\left[  y,x_{2n}\right]  =x_{2n}.$

\item $\tau_{2n+1}\left(  \lambda_{2}^{5},\cdots ,\lambda_{2}^{2n-1}\right)  $

$\hspace{1cm}\left[  x_{1},x_{k}\right]  =x_{k+1}, \ \ 2\leq k\leq2n-2,$

$\hspace{1cm}\left[  x_{k},x_{2n+1-k}\right]  =\left(  -1\right)  ^{k}x_{2n},
\ \ 2\leq k\leq n,$

$\hspace{1cm}\left[  y,x_{2+t}\right]  =x_{2+t}+\underset{k=2}{\overset
{\left[  \frac{2n-3-t}{2}\right]  }{\sum}}\lambda_{2}^{2k+1}x_{2k+1+t},
\ \ 0\leq t\leq2n-6,$

$\hspace{1cm}\left[  y,x_{2n-k}\right]  =x_{2n-k}, \ \ k=1,2,3,$

$\hspace{1cm}\left[  y,x_{2n}\right]  =2x_{2n}.$
\end{enumerate}
\end{proposition}

\subsubsection{Index of Lie algebras whose nilradical is $Q_{2n}$}

\begin{proposition}
Index of  $n$-dimensional Lie algebras whose nilradical is $Q_{2n}$ are
\begin{align*}
&\chi\left(  \tau_{2n+1}\left(  \lambda_{2}\right)  \right)     =1,\\
&\chi\left(  \tau_{2n+1}\left(  2-n,\varepsilon\right)  \right)     =1,\\
&\chi\left(  \tau_{2n+1}\left(  \lambda_{2}^{5},....,\lambda_{2}^{2n-1}\right)
\right)     =1.
\end{align*}

\end{proposition}

\begin{proof}
The corresponding matrix of the Lie algebra $\tau_{2n+1}\left(  \lambda
_{2}\right)  $ is of the form:

\begin{small}
$$\left(
\begin{tabular}
[c]{lllllll}%
$0$ & $x_{3}$ & $x_{4}$ & $\cdots$ & $0$ & $0$ & $-x_{1}$\\
$-x_{3}$ & $0$ & $0$ & $\cdots$ & $x_{2n}$ & $0$ & $-\lambda_{2}x_{2}$\\
$-x_{4}$ & $0$ & $0$ & $\cdots$ & $0$ & $0$ & $-\left(  n-\left(  n-1\right)
+\lambda_{2}\right)  x_{3}$\\
$\vdots$ & $\vdots$ & $\vdots$ & $\vdots$ & $\vdots$ & $\vdots$ & $\vdots$\\
$0$ & $-x_{2n}$ & $0$ & $\cdots$ & $0$ & $0$ & $-\left(  n-1+\lambda
_{2}\right)  x_{2n-1}$\\
$0$ & $0$ & $0$ & $\cdots$ & $0$ & $0$ & $-\left(  2n-3+2\lambda_{2}\right)
x_{2n}$\\
$x_{1}$ & $\lambda_{2}x_{2}$ & $\left(  n-\left(  n-1\right)  +\lambda
_{2}\right)  x_{3}$ & $\cdots$ & $\left(  n-1+\lambda_{2}\right)  x_{2n-1}$ &
$\left(  2n-3+2\lambda_{2}\right)  x_{2n}$ & $0$%
\end{tabular}
\ \right)  $$
\end{small}
Its rank is $2n$, then $\chi(\tau_{2n+1}\left(  \lambda_{2}\right)  )=1.$

The corresponding matrix of the algebra $\tau_{2n+1}\left(  2-n,\varepsilon
\right)  $ is of the form

$$\left(
\begin{tabular}
[c]{llllllll}%
$0$ & $x_{3}$ & $x_{4}$ & $\cdots$ & $x_{2n-1}$ & $0$ & $0$ & $-x_{1}%
-\varepsilon x_{2n}$\\
$-x_{3}$ & $0$ & $0$ & $\cdots$ & $0$ & $x_{2n}$ & $0$ & $-\left(  n-2\right)
x_{2}$\\
$-x_{4}$ & $0$ & $0$ & $\cdots$ & $-x_{2n}$ & $0$ & $0$ & $-\left(
n-2\right)  x_{3}$\\
$\vdots$ & $\vdots$ & $\vdots$ & $\vdots$ & $\vdots$ & $\vdots$ & $\vdots$ &
$\vdots$\\
$-x_{2n-1}$ & $0$ & $x_{2n}$ & $\cdots$ & $0$ & $0$ & $0$ & $-\left(
n-\left(  2n-1\right)  \right)  x_{2n-1}$\\
$0$ & $-x_{2n}$ & $0$ & $\cdots$ & $0$ & $0$ & $0$ & $-x_{2n}$\\
$0$ & $0$ & $0$ & $\cdots$ & $0$ & $0$ & $0$ & $0$\\
$x_{1}+\varepsilon x_{2n}$ & $\left(  n-2\right)  x_{2}$ & $\left(
n-2\right)  x_{3}$ & $\cdots$ & $\left(  n-\left(  2n-1\right)  \right)
x_{2n-1}$ & $x_{2n}$ & $0$ & $0$%
\end{tabular}
\ \right)  $$

\bigskip Its rank is $2n$, then $\chi(\tau_{2n+1}\left(  2-n,\varepsilon
\right)  )=1.$

Since the corresponding matrix of the algebra $\tau_{2n+1}\left(  \lambda
_{2}^{5},....,\lambda_{2}^{2n-1}\right)  $ is of rank $2n$ then the index is
$1$.
\end{proof}

\begin{remark}
The procedure described  in Remark \ref{regularvector} could be used  to compute the
regular vectors of  Lie algebras whose nilradical is $Q_{2n}$.
\end{remark}

\end{document}